\documentclass[onecolumn,draftcls]{IEEEtran}

\usepackage{url}
\usepackage{graphicx}
\usepackage{amsmath,amssymb} 
\usepackage{citesort}
\usepackage{color}
\usepackage{multirow}
\usepackage{url}

\usepackage{bm}
  \newcommand{\req}[1]    {(\ref{#1})}
  \newcommand{\refig}[1]  {Figure~\ref{#1}}
  \newcommand{\resec}[1]  {Section~\ref{#1}}
  \newcommand{\retab}[1]  {Table~\ref{#1}}

\begin{document}

\title{Walking dynamics are symmetric (enough)}

\author{M. Mert Ankaral{\i}, Shahin Sefati, Manu S. Madhav, Andrew Long, Amy J.
Bastian, and Noah J. Cowan \thanks{M. Mert Ankaral{\i}, Shahin Sefati, Manu S.
Madhav, and Noah J. Cowan are with the Dept. of Mechanical Eng., Johns Hopkins
University, Baltimore, MD, USA {\tt\small mertankarali@jhu.edu}}
\thanks{Andrew Long is with the Dept. of Biomedical Eng., Johns Hopkins
University, Baltimore, MD, USA } \thanks{Amy J. Bastian
    is with the Dept. of Neuroscience, Johns Hopkins University,
    Baltimore, MD, USA } }

\maketitle

\begin{abstract}

  Many biological phenomena such as locomotion, circadian cycles, and
  breathing are rhythmic in nature and can be modeled as rhythmic
  dynamical systems. Dynamical systems modeling often involves
  neglecting certain characteristics of a physical system as a
  modeling convenience. For example, human locomotion is frequently
  treated as symmetric about the sagittal plane.  In this work, we
  test this assumption by examining human walking dynamics around the
  steady-state (limit-cycle). Here we adapt statistical cross
  validation in order to examine whether there are statistically
  significant asymmetries, and even if so, test the consequences of
  assuming bilateral symmetry anyway. Indeed, we identify significant
  asymmetries in the dynamics of human walking, but nevertheless show
  that ignoring these asymmetries results in a more consistent and
  predictive model. In general, neglecting evident characteristics of
  a system can be more than a modeling convenience---it can produce a
  better model.

\end{abstract}

\section{INTRODUCTION}

The concept of symmetry has helped shape our understanding of
engineering and biology alike. The Roman text ``De architechura'' by
Vitruvius and the eponymous ``Vitruvian Man'' by Leornado Da Vinci
exemplify the influence of symmetry in animals and humans on man-made
works of art and engineering. Symmetry serves to simplify and reduce
model complexity, making it a powerful tool in computational and analytical
applications. The ubiquity of bilateral (left--right, sagittal plane)
symmetry in animals is genetically encoded~\cite{finnertyorigins2004},
and from an engineering point of view, building machines with
bilateral symmetry is justified by the fact that the left--right axis
is unbiased either by gravity or by direction of movement.  However,
genetic encoding of symmetry manifests itself imperfectly; numerous
factors, such as differences in contralateral limb lengths, dominance
of ``leggedness'' and handedness, and developmental processes break
perfect symmetry and enhance asymmetry.

Various measures and indices of asymmetry have been used to argue that
human locomotion is bilaterally symmetric or asymmetric (for reviews,
see \cite{sadeghisymmetry2000,hsiaoreview2010}). Symmetry is thought
to confer some advantages on motor abilities (e.g.\ improved energetic
efficiency
\cite{finleylearning2013,laigait2001,matteswalking2000}). The common
trend among previous work is the comparison of kinetic and/or
kinematic gait parameters between the right and left halves of the
body, i.e. joint
angles~\cite{hannahkinematic1984,forczekevaluation2012,karamanidissymmetry2003,reismanlocomotor2007},
ranges~\cite{stefanyshynright1994,gundersenbilateral1989} and
velocities~\cite{lawmicrocomputer1987}, stride
lengths~\cite{choderafootprint1973,choderaanalysis1974,gundersenbilateral1989,reismanlocomotor2007},
ground reaction
forces~\cite{hamillground1984,menardcomparative1992,vanderstraatensymmetry1978,herzogasymmetries1989},
EMG
profiles~\cite{carlsookinetic1973,arsenaultthere1986,marksanalysis1958,ounpuubilateral1989},
limbs forces and
moments~\cite{damholtasymmetry1978,balakrishanintegral1982,vaughanjoint1996,lathroplambachevidence2014},
or center-of-mass
oscillations~\cite{crowecharacterization1995,crowecharacterization1993,giakastime1997}.
However, as Sadeghi et.  al.~\cite{sadeghisymmetry2000} state, ``\dots
can we argue that it is acceptable to conclude that able-bodied gait
is asymmetrical just because of the existence of statistically
significant differences between two corresponding parameters (which we
call local asymmetry) calculated from the right and left limbs?''
During human walking, do steps from left-to-right and right-to-left
recover significantly differently from perturbations?  After all,
there are differences in leg dominance---e.g.\ preferred kicking
leg---that might lead to different responses from step-to-step.

Aside from demonstrating asymmetry (or not) in gait parameters, we
found no studies examining the potential benefits of \emph{neglecting}
evident asymmetries.  If there is a step-to-step dynamical
asymmetry, does fitting a model from stride to stride (two step)
rather than step to step (one step) better capture the dynamics of
human walking? Of course, no single physical system has perfect
symmetry. Thus symmetric models are inherently wrong for any
physical system, but may nevertheless be useful for simplifying both
the modeling and analysis.

``Essentially all models are wrong, but some are useful'' wrote George
E. P. Box in his seminal book \cite{boxempirical1987}. According to
Box, the important practical concern regarding the models of physical
phenomena is ``how wrong do they have to be to not be useful?'' With
regard to bilateral asymmetry in human walking, we attempt to frame
this concern as follows. How wrong is it to neglect asymmetry from a
statistical point of view?  And how useful is symmetric modeling in
terms of predictive power and simplicity?  In most cases, correctness
and usefulness are directly related, and they are tested
simultaneously. However, in the context of data-driven modeling of
human walking dynamics, the ``wrongness'' and ``usefulness'' of assuming
symmetry are related but have critical, nuanced differences. The
methods presented in this paper allows us to independently
(statistically) address these differences.

In this paper, we test the assumption of bilateral symmetry in the
dynamics of human walking. As an example, consider fitting linear
models to two distinct data sets (e.g. ``left steps'' and ``right
steps'') and testing these models in terms of their respective ability
to predict isolated validation data from just the one of the data
sets, say ``left steps''. If walking were perfectly symmetric, both
the left-step (``correct'') model and right-step (``wrong'') model
would perform indistinguishably in left-step validation. However, we
show that there are statistically significant asymmetries in the
dynamics of human walking in healthy subjects in the sense that the
``wrong'' model performs statistically worse than the ``correct''
model in validation. Despite these asymmetries, we also show that a
more consistent and predictive model of the dynamics is obtained by
assuming symmetry, and pooling all the data from both left and right
steps to form a generic model. Quite surprisingly, this fit
significantly out-performs the mapping fitted to only left steps even
when predicting left-step data.  This is good news because, in addition
to our finding that it is statistically better to neglect asymmetry,
it is also practically and theoretically convenient to assume
symmetry. These advantages lead us to conclude that the assumption of
symmetry in walking dynamics, though clearly wrong in a platonic
sense, is nevertheless more useful for all practical purposes.

\subsection{Modeling the Rhythmic Dynamics}
\label{sec:rhythm}

Our approach to analyzing and modeling walking involves treating the
underlying behavior as a finite-dimensional nonlinear rhythmic
dynamical system operating around a stable limit cycle. This type of
modeling approach has been successful for robotic
\cite{altendorferstability2004a,chevallereauasymptotically2009,ankarali_saranli.chaos2010,depenn2015}
and biological systems
\cite{hurmuzlumeasurement1994,seyfarthswing-leg2003,holmesdynamics2006,revzenfinding2012,ankaralihaptic2014}. 
A limit cycle is an isolated periodic trajectory that is a solution to
the equations governing the dynamical system
\cite{guckenheimernonlinear1991}. A limit cycle is said to be stable
if all trajectories in a sufficiently small neighborhood of the limit
cycle converge to it.

We further use Poincar\'e theory in our analysis of rhythmic walking
dynamics. A Poincar\'e return map
\cite{guckenheimernonlinear1991,holmesdynamics2006} is a mapping from
a transverse section $S$ back to itself, obtained by tracing the
consecutive intersections of the state trajectories with the section
$S$.  This return map reduces the continuous rhythmic dynamical system
to a nonlinear discrete dynamical system that preserves many
properties of the behavior. The specific Poincar\'e section that we
adopt for human walking is the heel strike event, as explained in
\resec{sec:poincare}.

The intersection of the limit cycle with the Poincar\'e section is an
isolated fixed point of the return map. The limit cycle is
asymptotically stable if and only if this fixed point is stable. Our
second modeling approximation is based on the Hartman–-Grobman theorem
(or linearization theorem), which states that local flow around any
hyperbolic fixed point is homeomorphic to the one governed by its
linearization around the fixed point itself. Thus as detailed in
\resec{sec:poincare}, we fit linear models to walking trajectories on
the Poincar\'e maps.

\subsection{Limit-Cycle Dynamics and Symmetry}
Here, we define symmetry in the context of limit-cycle modeling of
walking and consider what kind of symmetries (and asymmetries) can be
addressed using this approach.

In our modeling approach, there are two core elements: the limit-cycle of the
rhythmic system, which characterizes the steady-state behavior, and the dynamics
(both deterministic and stochastic) around the limit-cycle. In this paper, we are
interested in the latter. 

Beyond its utility for approximation, bilateral symmetry of the
(steady-state) limit-cycle trajectory may have physiological
significance, such as reducing metabolic
cost~\cite{matteswalking2000,finleylearning2013,laigait2001}.  Indeed,
the kinematic and kinetic variables that are the focus of the majority
of studies that address human~\cite{sadeghisymmetry2000} or
animal~\cite{muircomplete1999,pourcelotkinematic1997} locomotor
symmetry are steady-state (periodic) variables that correspond to the
limit cycle of a dynamic model.

Here, we consider the dynamics near, but \emph{off of} the limit
cycle, using data from Poincar\'e sections to estimate return
maps. Hence our analysis is not based on the steady-state parameters
of gait, but how the gait deviates from and recovers to these
steady-state parameters. To the best of our knowledge, this is the
first study that analyzes the \emph{dynamical} symmetry of biological
rhythmic systems. We validate our methods in data from normal human
walking experiments. These methods are also applicable to robotic or
biological locomotor behavior with approximately symmetric gait
patterns.

\section{METHODS}

The system of interest is human treadmill walking. This data set is obtained
for eight healthy young adult participants, at three different belt speeds
(0.5, 1.0, and 1.5 m/s).  We required the subjects to cross their arms  in
order to continuously record the marker positions. The Johns Hopkins
Institutional Review Board approved all protocols and all subjects gave
informed written consent prior to participation.

\subsection{Kinematic Data}

\begin{figure}[ht] \centering
  \includegraphics{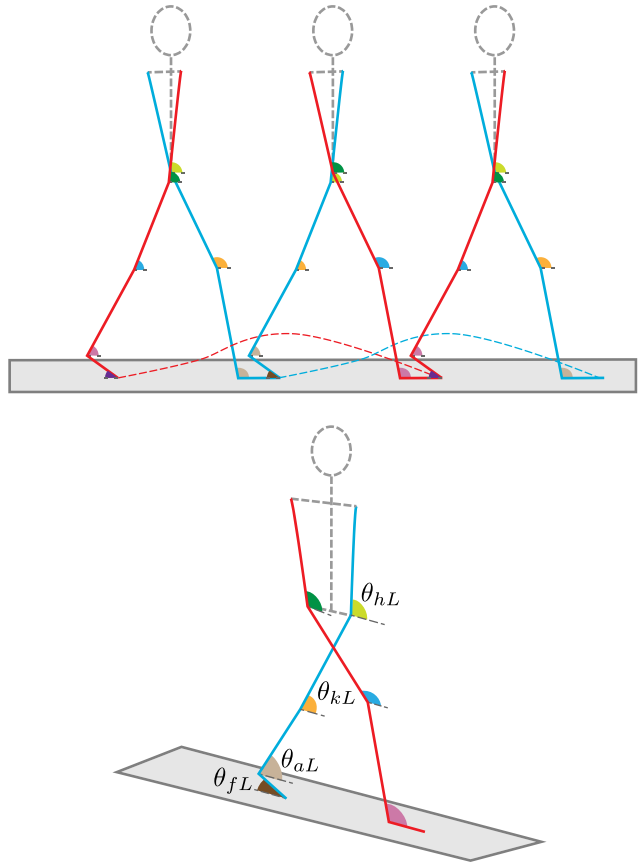}
  \caption{Visualization of leg angle vector $\bm{\theta}$. The left
    leg (blue) and right leg (red) alternate between stance and swing
    phase over the course of a stride. The variables $\theta_{fL}$,
    $\theta_{aL} $, $\theta_{kL}$, and $\theta_{hL}$ correspond to the
    left foot, ankle, knee, and hip angles respectively. The
    corresponding right leg angles, $\theta_{fR}$, $\theta_{aR} $,
    $\theta_{kR}$, and $\theta_{hR}$, are not labeled. The $8$ leg
    angles and their respective angular velocities form the
    $16$-dimensional state vector.}
  \label{fig:walkingState}
\end{figure}

We placed infrared (IR) markers on subjects' left and right shoulder, hip,
knee, ankle, and toe. Markers were tracked in 3D using Optotrak (Northern
Digital) at 100 Hz. The marker data were used to calculate the four sagittal
plane angles on each side as illustrated in Figure~\ref{fig:walkingState}. The
raw angular data were smoothed with a fifth-order, zero-phase-lag (non-causal)
Butterworth filter ($5^{th}$ order with a cut-off frequency of $10 \mathrm{Hz}$) 
to remove measurement noise and ease angular velocity
estimation. In order to estimate angular velocities, a central difference
filter and another zero-phase-lag Butterworth filter was applied to the
smoothed angular data similar to the methods adopted in the biomechanics 
literature \cite{millercontinuous2008,kurzdifferences2012}. We assume that 
the smoothed angles (8) and angular velocities (8) form a 16 dimensional 
state space for walking. The state vector includes angles (rad),
\begin{equation}
  \label{eq:1}
  \begin{aligned} 
    \bm{\theta_L}(t) &= \begin{bmatrix} 
      \theta_{fL} 
      & \theta_{aL} 
      & \theta_{kL} 
      & \theta_{hL} \end{bmatrix}^T, \\
    \bm{\theta_R}(t) &= \begin{bmatrix}
      \theta_{fR} 
      & \theta_{aR} 
      & \theta_{kR} 
      & \theta_{hR} 
    \end{bmatrix}^T, \\
    \bm{\theta}(t) &= \begin{bmatrix}
      \bm{\theta_L}(t) \\ \bm{\theta_R}(t)
    \end{bmatrix}
  \end{aligned}
\end{equation}
and angular velocities (rad/s),
\begin{equation}
  \begin{aligned} 
    \dot{\bm{\theta}}_{\bm{L}}(t) &= 
    \begin{bmatrix} 
      \dot{\theta}_{fL} 
      & \dot{\theta}_{aL} 
      & \dot{\theta}_{kL} 
      & \dot{\theta}_{hL}
    \end{bmatrix}^T, \\
    \dot{\bm{\theta}}_{\bm{R}} (t) &=
    \begin{bmatrix}   \dot{\theta}_{fR} 
      & \dot{\theta}_{aR} 
      & \dot{\theta}_{kR} 
      & \dot{\theta}_{hR}  
    \end{bmatrix}^T, \\
    \dot{\bm{\theta}}(t) &= 
    \begin{bmatrix}
      \dot{\bm{\theta}}_{\bm{L}}(t) \\ \dot{\bm{\theta}}_{\bm{R}}(t)
    \end{bmatrix}.
  \end{aligned}
\end{equation}
Subscripts $f$, $a$, $k$, and $h$ stand for foot, ankle, knee, and
hip, respectively. $L$ and $R$ mnemonically denote the left and right
legs.

In order to analyze the data independently from the physical units,
the state space was non-dimensionalized based on the time constant
associated with pendular walking
\cite{donelanmechanical2002,ankarali_saranli_AR2011,saranli_arslan_ankarali_morgul.nd2010,blickhansimilarity1993}:
\begin{equation}
  \label{eq:4}
  \begin{aligned} 
    \bar{\bm{\theta}} &= \bm{\theta} \\
    \bar{\dot{\bm{\theta}}} &= \dot{\bm{\theta}} \sqrt{\frac{l_0}{g}},
  \end{aligned}
\end{equation}
where the bar represents the corresponding non-dimensionalized
variable, $g$ is the gravitational acceleration, and $l_0$ is the leg
length of the subject, which is estimated from the marker
positions on right hip and ankle.

\subsection{Events and Section Data}
\label{sec:poincare}
The treadmill used in this study features a split belt\footnote{While
  the belts of the treadmill can be driven at different speeds, this
  study addresses bilateral symmetry, so both belts were driven at the
  same speed.} that mechanically decouples the vertical ground
reaction forces caused by each foot. Each belt is instrumented with a
separate load cell, facilitating the estimation of the timing of
heel-strike events. We chose heel-strike events as Poincar\'e sections
for the analyses.

Let $t[k]$ be the detected times of heel-strike events, where $k \in
\mathcal{K}=\{1,2,3,\dots,k_{\mathrm{max}}\}$, with $k_{\mathrm{max}}$
being the total number of heel-strike events of both legs in one
walking trial. For example, if the first heel-strike event ($k=1$)
corresponds to the left leg, sets of odd ($\mathcal{K}_L$) and even
($\mathcal{K}_R$) integer indices from 1 to $k_{\mathrm{max}}$
correspond to the left and right heel-strike events, respectively,
such that $\mathcal{K}=\mathcal{K}_L \cup \mathcal{K}_R$. Over one
stride of walking, there are two Poincar\'e sections of interest at
heel-strike events. The measurement of the state vector at these
Poincar\'e sections is given as follows:
\begin{equation}
  \textbf{z}[k] =  
  \begin{bmatrix}
    \bar{\bm{\theta}}(t[k]) \\
    \bar{\dot{\bm{\theta}}}(t[k]) 
  \end{bmatrix} .
\end{equation}
During steady-state walking and in the absence of noise, the periodic
orbit would remain on the limit cycle:
\begin{equation}
  \begin{aligned}
    \bm{z}[m] &=\bm{\mu}_L\ , \ \forall m \in \mathcal{K}_L  
    \\ 
    \bm{z}[m'] & =\bm{\mu}_R \ , \ \forall m' \in \mathcal{K}_R    
  \end{aligned}
\end{equation}
where $\bm{\mu}_L$ and $\bm{\mu}_R$ are the fixed-points with respect
to each of the two distinct Poincar\'e sections. Note that assuming
bilateral symmetry implies that these two fixed points are identical
up to a relabeling
\cite{westervelthybrid2003,chevallereauasymptotically2009,duindammodeling2009,leetemplates2008}. This
relabeling can be expressed as a linear mapping of right-heel-strike
coordinates:
\begin{equation}
  M:\begin{bmatrix} \bar{\bm{\theta}}_L(t[k]) 
      \\ \bar{\bm{\theta}}_R(t[k]) 
      \\ \bar{\dot{\bm{\theta}}}_L(t[k]) 
      \\ \bar{\dot{\bm{\theta}}}_R(t[k])  \end{bmatrix}
  \mapsto
   \begin{bmatrix} \bar{\bm{\theta}}_R(t[k]) 
      \\ \bar{\bm{\theta}}_L(t[k]) 
      \\ \bar{\dot{\bm{\theta}}}_R(t[k]) 
      \\ \bar{\dot{\bm{\theta}}}_L(t[k])  \end{bmatrix}, \ \forall k \in \mathcal{K}_R ,
\end{equation}
where
\begin{equation}
  M = \begin{bmatrix}
    0 & I_{4\times4} & 0 & 0 \\ 
    I_{4\times4} & 0 & 0 & 0 \\
    0 & 0 & 0 & I_{4\times4} \\ 
    0 & 0 & I_{4\times4} & 0 \\
  \end{bmatrix} .
\end{equation}

\begin{figure*}[tb!] \centering
  \includegraphics{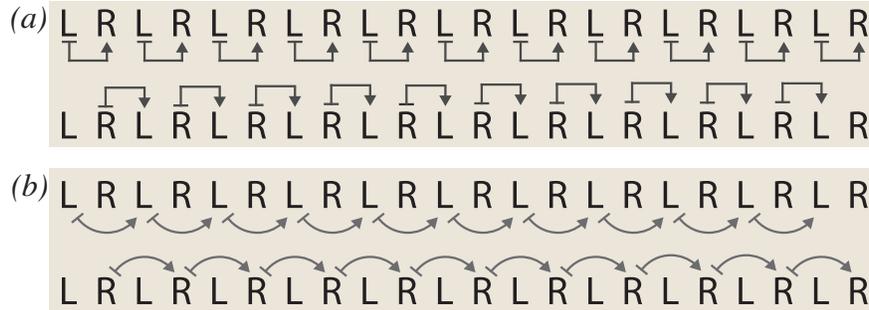}
  \caption{Types of input--output pairs analyzed in this paper. $L$
    and $R$ represent the Poincar\'e sections associated with
    heel-strike events of the left and right legs.  (a) Left-to-right
    step maps (top) and right-to-left step maps (bottom). Step maps
    are denoted using straight arrows. (b) Left-to-left stride maps
    (top) and right-to-right stride maps (bottom). Stride maps and
    step maps are distinguished throughout the paper by shape
    (straight versus curved arrows, respectively). }
  \label{fig:data}
\end{figure*}

As explained in \resec{sec:rhythm} our approach to modeling human
walking centers around fitting linear maps between Poincar\'e sections
around the associated fixed points. First, we estimated the
fixed points via
\begin{equation}
  \label{eq:5}
  \begin{aligned}
    \hat{\bm{\mu}}_L & = \frac{1}{|\mathcal{K}_L|} \sum\limits_{k\in\mathcal{K}_L} \bm{z}[k], \\
    \hat{\bm{\mu}}_R & = \frac{1}{|\mathcal{K}_R|} \sum\limits_{k\in \mathcal{K}_R} \bm{z}[k],
  \end{aligned}
\end{equation}
where $|\mathcal{K}_L|$ and $|\mathcal{K}_R|$ denote the cardinality
of sets $\mathcal{K}_L$ and $\mathcal{K}_R$ respectively. Note that
\emph{kinematic} asymmetry could be measured directly in terms of the
difference between respective fixed points $\hat{\bm{\mu}}_L$ and $
\hat{\bm{\mu}}_L$. While potentially of interest, the current paper
focuses on \emph{dynamical} asymmetry (measured in terms of the section
maps), and thus we computed the residuals by subtracting the estimated
fixed points from the section data:
\begin{equation}
  \begin{aligned}
    \bm{q}_L[k] & =  \bm{z}[k] - \hat{\bm{\mu}}_L, & k\in \mathcal{K}_L,\\
    \bm{q}_R[k] & = \bm{z}[k] - \hat{\bm{\mu}}_R, & k\in \mathcal{K}_R.
  \end{aligned}
\end{equation}
Section maps were estimated using these residuals. A section map from
$\bm{q}_L$ to the subsequent $\bm{q}_R$ is denoted as $L \mapsto R$. We
fit two categories of section maps: step-to-step ($L \mapsto R$ and $R
\mapsto L$) and stride-to-stride ($L \mapsto L$ and $R \mapsto
R$). These sections maps are illustrated in \refig{fig:data}(a). 

\subsection{Fitting Section Maps}
To fit the section maps for each category explained above, we stack
all the appropriate residuals ($\bm{q}_L$ and/or $\bm{q}_R$) in
matrices $X$ (input) and $Y$ (output):
\begin{equation} X = [x_1 , \cdots , x_N]^T , \ Y = [y_1 , \cdots ,
  y_N]^T ,\ x_i , y_i \in \mathbb{R}^d,
  \label{eq:XY}
\end{equation}
where $x_i$ and $y_i$ represent residuals ($\bm{q}$) from sections
evaluated in the data. For example, to fit the $L\mapsto R$
step-to-step map, one would set the columns of $X$ and $Y$ as follows:
\begin{equation}
  \label{eq:6}
  \begin{aligned}
    x_1 &= q_L[1], &y_1 &= q_R[2],\\
    x_2 &=q_L[3], &y_2 &= q_R[4],\\
    &\,\vdots & &\,\vdots\\
    x_N &=q_L[2N - 1], &y_N &= q_R[2 N].
  \end{aligned}
\end{equation}
The linear section maps are modeled as $y_i = A x_i + \delta_i , \ \forall i$,
%% 
%\begin{align}
%  y_i &= A x_i + \delta_i \ , \ \forall i ,\\
%  Y &= X A^T + \Delta ,
%\end{align}
where $\delta_i$ is additive noise.
The section map can be estimated via least squares:
\begin{align}
  \label{eq:returnMapFitting}
  \hat{A} = \left( X^\dagger Y \right)^T.
\end{align} 
where $X^\dagger$ is the Moore-–Penrose pseudoinverse
of $X$.

\section{STATISTICAL APPROACH}

Here, we tailor Monte-Carlo cross validation to examine symmetry in
walking dynamics. (A similar technique based on bootstrap sampling
produces qualitatively similar
results~\cite{ankaralivariability2015}.)

\subsection{Test of Symmetry Using Monte Carlo Cross-Validation}
\label{sec:montecarlo}

Classical cross validation (CV) involves fitting a model to a
\textit{training set} of input--output data and validating the model
by comparing its predictions on a complementary \textit{test set} of
input--output data. In classical CV, there are $n$ pairs of
input--output data which are then split into a training (fitting) set
($n_f$ pairs) and complementary test (validation) set ($n_v = n - n_f$
pairs). The training set is used for model fitting. The fitted model
is then applied to the inputs of the test set to generate output
predictions; the error metric between the predicted and actual outputs
is the \textit{cross-validation error} (CVE). The CVE is used to
evaluate the performance of the model. CV methods are commonly used
for selecting models based on their predictive
ability~\cite{shaolinear1993,madhavclosed-loop2013,raolinear2005,yangconsistency2007}. A
review by Arlot and Celisse~\cite{arlotsurvey2010} summarizes
different cross-validation methods and discusses their advantages and
limitations.

%A critical question when using a CV method is how to split the
%data~\cite{shaolinear1993}. Assuming no replacement, there exist
%$\binom{n}{n_v}$ different ways of splitting the data set.  The most popular
%CV method, often called leave-one-out cross validation (LOOCV), uses $n_v =
%1$, because it incurs the least computational cost. However, LOOCV is
%inconsistent (asymptotically biased)~\cite{shaolinear1993}, and its
%performance is poor in practice if the sample size is large.  However for $n_v
%\gg 1$ and $(n-n_v) \gg 1$, computing cross validation for all
%$\binom{n}{n_v}$ possibilities can be computationally expensive. 

%One way of
%decreasing the computational complexity is to apply $k$-fold cross validation
%\cite{kohavstudy1995,arlotsurvey2010}.
%
%In $k$-fold CV, the data set is initially split into $k$ mutually exclusive,
%equal-sized subsets, and a LOO-type analysis is performed at the level of the
%subsets. Note that $k$-fold CV is computationally much more feasible than
%testing all possible ways of splitting up the data, but there is a trade-off:
%picking small $k$ increases the variance of the CVE estimate, while picking
%large $k$, which tends toward LOOCV, is biased for  large sample sizes.
%
%A common way to overcome the artificial trade-off between bias and
%variance imposed by $k$-fold CV, and still offer computational
%tractability, 

The method that we present in this paper is based on Monte-Carlo cross
validation (MCCV)~\cite{shaolinear1993}.  MCCV randomly splits the
data $m$ times with fixed $n_f$ and $n_v$ (size of training and test
sets respectively) over the $m$ iterations. For each iteration, the
CVE is computed using the respective training and test sets; the
overall CVE is estimated using the mean of these $m$ CVEs.
%If $ \frac{n}{n_v} \ll m \ll \binom{n}{n_v}$, MCCV estimates the CVE with
%a computationally feasible sample size.
%, and with negligible variance compared to the estimate provided by $k$-fold CV. 
%Also, since we fit the model more times ($m \gg k$), MCCV provides a better mechanism for
%estimating model uncertainty than $k$-fold CV.

As mentioned before, the model being fit to input--output data in our
case is a linear map. Suppose there are $n$ pairs of input--output
data, $(x_i,y_i)$, where $i\in \mathcal{I} = \{1,2,3,\dots,n
\}$. Split this data into a training set $\mathcal{F}\ni(x_f,y_f)$
comprising $n_f$ pairs, and test set $\mathcal{V}\ni(x_v,y_v)$
comprising $n_v$ pairs.  Define $X_{\mathcal{F}}$ as the matrix whose
rows are $x_f^T$, and define $Y_{\mathcal{F}}$, $X_{\mathcal{V}}$, and
$Y_{\mathcal{V}}$ similarly. We use the following definition of CVE
from $\mathcal{F}$ to $\mathcal{V}$:
\begin{align} 
  \mathrm{CVE}_{\mathcal{F} \mapsto \mathcal{V}} := 
  \frac{||Y_{\mathcal{V}}
    - X_{\mathcal{V}} A_{\mathcal{F}}^T||^2}{||Y_{\mathcal{V}}||^2}
  \label{eq:CVE}
\end{align}
where $|| \cdot ||$ denotes the Frobenius norm and
\begin{align}
  A_\mathcal{F} := \left( X_{\mathcal{F}}^\dagger Y_{\mathcal{F}}
  \right)^T
  \label{eq:model}
\end{align}
is the least-squares solution \eqref{eq:returnMapFitting} given 
the training data $\mathcal{F}$.

%\subsection{Extending Monte Carlo Cross-Validation}
%\label{sec:extended}

We tailor classical Monte Carlo CV for systems that may exhibit
discrete symmetry. We focus our discussion and notation on human
walking, but these methods are applicable to other forms of locomotion
that involve nearly bilaterally symmetric gaits, e.g. walking and
trotting \cite{blickhansimilarity1993}, but not clearly asymmetric gaits, e.g. galloping
\cite{collinscoupled1993}. In classical Monte Carlo CV, at each
iteration, one CVE is computed using~\req{eq:CVE}, where as in our
CV method we compute three types of CVEs. Each CVE
computation uses the same test set, but the models are fit using three
different training sets.

Each application of extended CV method requires a ``normal'' set,
$\mathcal{N}$, and an equal size ``mirrored'' set, $\mathcal{M}$. In
this paper, we analyze four different $(\mathcal{N},\mathcal{M})$
pairs, which are generated using the input--output data types
illustrated in~\refig{fig:data}, i.e.\ step-to-step transitions
($\lbrace L \mapsto R \rbrace$ and $\lbrace R \mapsto L \rbrace$) and
stride-to-stride transitions ($\lbrace L \mapsto L \rbrace$ and
$\lbrace R \mapsto R \rbrace$).  For example, if the normal data set
comprises the left-to-right step transitions,
$\mathcal{N}=\lbrace L \mapsto R \rbrace$, the associated mirrored
dataset is $\mathcal{M}=\lbrace R \mapsto L \rbrace$.  Similarly, for
strides, if $\mathcal{N} =\lbrace L \mapsto L \rbrace$ represents the
set of all transitions from left heel strike to the subsequent left
heel strike, then $\mathcal{M}=\lbrace R \mapsto R \rbrace$ are the
corresponding right-to-right transitions.  All
$(\mathcal{N},\mathcal{M})$ combinations are listed in \retab{tab:cv}.

The normal and mirrored sets \emph{each} include $n$ mutually
exclusive input--output pairs, denoted by $ (x_i , y_i) \in
\mathcal{N}$ and $(\hat{x}_i,\hat{y}_i)\in\mathcal{M}$, respectively
where $i\in\mathcal{I} = \lbrace 1 , 2 , \cdots, n \rbrace$.  Each
iteration of extended CV randomly splits this index set
$\mathcal{I}$ into a training index set $\mathcal{I}_f$ and test index
set $\mathcal{I}_v$ in a manner identical to classical CV:
$\mathcal{I}_f \cup \mathcal{I}_v = \mathcal{I}$ and $\mathcal{I}_f
\cap \mathcal{I}_v =\emptyset$.  The three types of CVE computations
described below draw the test set from the normal dataset:
\begin{align}
  \mathcal{V} &= \{(x_v,y_v) \in\mathcal{N} \, | \,  v \in \mathcal{I}_v \} .
\end{align}

\textbf{\emph{Normal Cross Validation (NCV)}} is the same as classical
Monte Carlo CV, in that the training data are also drawn from
$\mathcal{N}$:
\begin{align}
  \mathcal{F}_{\mathrm{NCV}} &= \{(x_f,y_f) \in\mathcal{N} \, | \,  f \in \mathcal{I}_f \} .
\end{align}
This set is used for fitting the linear model
$A_{\mathcal{F}_{\mathrm{NCV}}}$ using \req{eq:model}.  Given
$A_{\mathcal{F}_{\mathrm{NCV}}}$, the CVE is computed on the common
test set $\mathcal{V}$ using \req{eq:CVE}. For NCV, the mirrored data set
$\mathcal{M}$ is not used.

\textbf{\emph{Mirrored Cross Validation (MCV)}} draws the training
data from the mirrored dataset $\mathcal{M}$, using the same training
index set $\mathcal{I}_f$ as NCV:
\begin{align}
  \mathcal{F}_{\mathrm{MCV}} &= \{(\hat{x}_f,\hat{y}_f) \in\mathcal{M}\, | \,  f \in \mathcal{I}_f \} .
\end{align}
As before, this set is used for computing the linear model
$A_{\mathcal{F}_{\mathrm{MCV}}}$ using \req{eq:model}. The common test set
$\mathcal{V}$ is used for computing the CVE.  In MCV, we are using the
``wrong'' training data (mirrored), which will be critical to detect dynamical
asymmetry in walking. Note that the size (and in fact the indices) of the
training data in both MCV and NCV are the same.

\textbf{\emph{Combined Cross Validation} (CCV)} uses training data
that is the union of the training sets from NCV and MCV:
\begin{align}
  \mathcal{F}_{\mathrm{CCV}} &= \mathcal{F}_{\mathrm{NCV}} \cup \mathcal{F}_{\mathrm{MCV}}.
\end{align}
And, as before, this data subset is used to fit the linear model
$A_{\mathcal{F}_\mathrm{CCV}}$ and the common test set $\mathcal{V}$ is used
for calculating the CVE.  Thus the model is fitted on data pooled from both
$\mathcal{N}$ and $\mathcal{M}$, while the test data remain the same. Note
that CCV uses twice as much data for fitting as either NCV or MCV.

\begin{figure}[ht] \centering
  \includegraphics{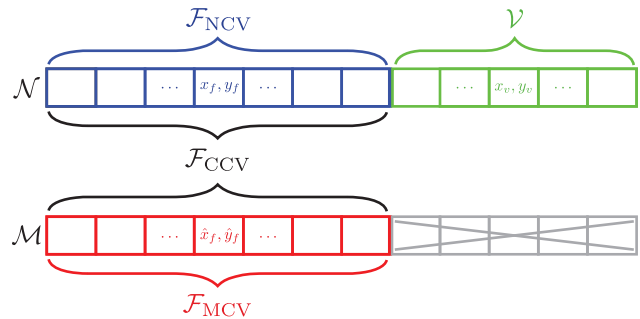}
  \caption{Illustration of the subsets $\mathcal{V}$,
    $\mathcal{F}_{\mathrm{NCV}}$, $\mathcal{F}_{\mathrm{MCV}}$, and
    $\mathcal{F}_{\mathrm{CCV}}$ after random splitting during an
    iteration of extended cross validation methods. The normal
    dataset, $\mathcal{N}$, is randomly split into the normal training
    set $\mathcal{F}_{\mathrm{NCV}}$ and the common test set
    $\mathcal{V}$.  $\mathcal{F}_{\mathrm{MCV}}$ shares the same
    indices as $\mathcal{F}_{\mathrm{NCV}}$ but is drawn from the
    mirrored dataset, $\mathcal{M}$. The training set for the CCV is
    simply the union of the other two training sets:
    $\mathcal{F}_{\mathrm{CCV}} = \mathcal{F}_{\mathrm{NCV}} \cup
    \mathcal{F}_{\mathrm{MCV}}$.  Note that the subset $\mathcal{M}
    \setminus \mathcal{F}_{\mathrm{MCV}}$ (greyed out) is not used in
    any of the three CV computations.  }
  \label{fig:method}
\end{figure}

\begin{figure*}[ht] \centering
  \includegraphics{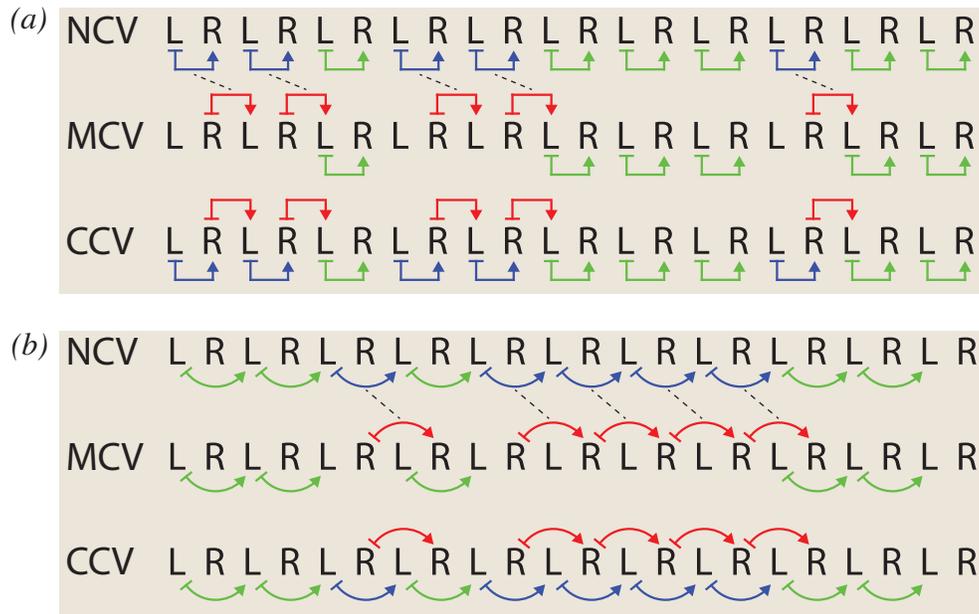}
  \caption{Illustration of extended CV dataset partitioning. (a) For
    step-to-step data, the normal dataset ($\mathcal{N}$) comprises
    all left-to-right step ordered pairs, whereas the mirrored dataset
    ($\mathcal{M}$) comprises all right-to-left step ordered
    pairs. (b) For stride-to-stride data, the normal dataset
    ($\mathcal{N}$) comprises all left-to-left stride ordered pairs,
    whereas the mirrored dataset ($\mathcal{M}$) comprises all
    right-to-right stride ordered pairs.  In both cases, for each
    iteration, a common test set ($\mathcal{V}$, green arrows), used
    for all CV methods, is randomly sampled from the normal
    dataset. The training sets, however, are unique to each method.
    NCV: the remainder of the normal dataset is used for training
    ($\mathcal{F}_{\mathrm{NCV}}$, blue arrows).  MCV: the training
    set ($\mathcal{F}_{\mathrm{MCV}}$, red arrows) is obtained using
    the same indices (dashed lines) as for
    $\mathcal{F}_{\mathrm{NCV}}$. CCV: the union of the test sets for
    NCV and MCV, comprise the combined training data
    ($\mathcal{F}_{\mathrm{CCV}} = \mathcal{F}_{\mathrm{NCV}} \cup
    \mathcal{F}_{\mathrm{MCV}}$, red and blue arrows).}
  \label{fig:extendedcrossstep}
\end{figure*}

\refig{fig:method} illustrates the set partitioning for one iteration of
extended cross validation for a general dataset.
\refig{fig:extendedcrossstep}(a) illustrates one iteration of the extended
cross validation algorithms on step data, where the normal dataset is
$\mathcal{N} = \lbrace L \mapsto R \rbrace$ and the mirrored dataset is
$\mathcal{M} = \lbrace R \mapsto L \rbrace$. \refig{fig:extendedcrossstep}(b)
illustrates one iteration of the extended cross validation algorithms on
stride data where $\mathcal{N} = \lbrace L \mapsto L \rbrace$ and $\mathcal{M}
= \lbrace R \mapsto R \rbrace$.	

The comparison between NCV and MCV will be critical for statistically
testing the symmetry of human walking. Since both NCV and MCV have
training sets of the same size and MCV uses mirrored data, the
difference in CVEs offer a direct measure of dynamical asymmetry. 
For a symmetric system, NCV and MCV errors should be
statistically indistinguishable. If there are asymmetries, we should
observe higher MCV errors than NCV errors.

However, this comparison alone is not enough to address all concerns
because the main advantage of assuming symmetry is that we double the
amount of data by combining the normal and mirrored data sets and
fitting a single model. We introduce a potential bias by neglecting
the asymmetries in the behavior, but reduce the variance in the
estimation by simply doubling the amount of data used for fitting.
From this perspective, comparison between NCV and CCV will be critical
for statistically testing predictive powers of asymmetric and
symmetric modeling approaches, which is an effective way testing the
``usefulness'' of the symmetry assumption.

\section{RESULTS}
\label{sec:results}

The results presented here are based on the methods presented in
\resec{sec:montecarlo} which rely on Monte-Carlo sampling and
cross-validation. The results presented below were qualitatively
similar (and stronger in one case) to those obtained using the
bootstrap method presented in \cite{ankaralivariability2015}.

We set the sample size of Monte-Carlo iterations to $m=1000$ based on
pilot experiments which showed that increasing the sample size beyond
this had a negligible effect on cross validation error. In each
iteration, $20\%$ ($\frac{n_v}{n} = 0.2$) of the normal dataset,
$\mathcal{N}$, was withheld for validation. Training sets for the
three CV computations were drawn from the remaining data according to
the procedure detailed in \resec{sec:montecarlo}.

\subsection{Symmetric vs.\ Asymmetric Modeling}
\label{sec:symvsasym}

The question being addressed in this paper is not just the symmetry
versus asymmetry of the dynamics of human walking, but also the
statistical consequences of choosing one approach over the other. We
applied our cross-validation method (\resec{sec:montecarlo}) to
expose these consequences.

\subsubsection{Step Maps}
\label{sec:stepmaps}

To apply the extended CV to step-to-step transitions, we
analyzed both combinations of normal and mirrored data:
$\left(\mathcal{N},\mathcal{M} \right) = \left(\lbrace L \mapsto R
  \rbrace , \lbrace R \mapsto L \rbrace\right)$ and $\left(
  \mathcal{N},\mathcal{M} \right) = \left(\lbrace R \mapsto L \rbrace
  , \lbrace L \mapsto R \rbrace \right)$ (Table \ref{tab:cv}). For
each category of cross validation---NCV, MCV, and CCV---we averaged
the errors for both combinations of $\left(\mathcal{N},\mathcal{M}
\right)$.  

\begin{figure*}[ht] \centering
  \includegraphics[scale=1]{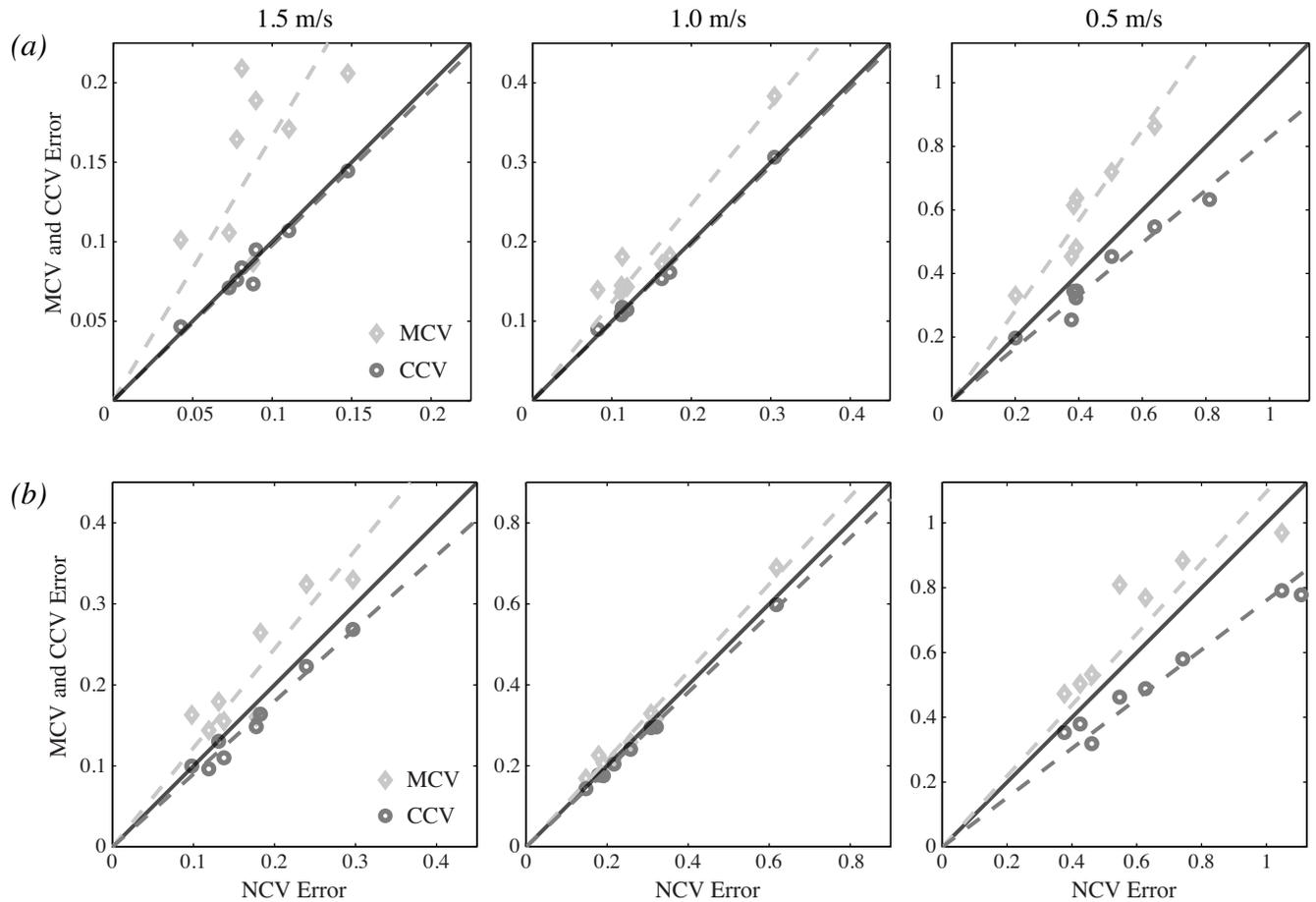}
  \caption{Walking is asymmetric, but neglecting this by training a
    model on the combined data can nevertheless improve fitting. (a)
    Step-to-step maps.  At all speeds, the mean mirrored cross
    validation errors (light grey diamonds) were significantly worse than
    for normal cross validation, indicating that steps were indeed
    asymmetric.  Despite this left--right asymmetry, the mean combined
    cross validation errors were not significantly different than for
    normal cross validation at the two fasted walking speeds tested,
    and, more surprisingly, were actually \emph{lower} at the slowest
    walking speed. The slopes of the fitted lines (dashed) determine
    the relative increase ($m>1$) or decrease ($m<1$) in CV error
    relative to the NCV error. (b) Stride-to-stride maps. By the same
    statistical measure, strides were also asymmetric at all speeds,
    but less substantially so. Moreover, the mean CCV error was lower
    than mean NCV error at all speeds.}
  \label{fig:CrossCrossBoth}
\end{figure*}

\refig{fig:CrossCrossBoth}(a) compares MCV and CCV errors to NCV error
from step-to-step data. MCV errors are (statistically) significantly
higher than NCV errors at all speeds ($p_{1.5 m/s} = 0.004$, $p_{1
  m/s} = 0.008$ and $p_{0.5 m/s} = 0.004$; one-sided Wilcoxon
rank-sign test).  This shows that our dataset is indeed dynamically
asymmetric between $L \mapsto R$ and $R \mapsto L$.

The comparison of CCV and NCV errors illuminates a
different perspective (Figure \ref{fig:CrossCrossBoth}(a)). For speeds of
$1.5$ m/s and $1.0$ m/s, CCV and NCV errors were statistically
indistinguishable ($p_{1.5 m/s} = 0.38$ and $p_{1 m/s} = 0.84$;
Wilcoxon rank-sign test), suggesting that for these speeds, the
predictive power of a model that assumes symmetry is just as great as
one that embraces the asymmetry. More surprisingly, the average CCV
error for the slowest speed tested was (statistically) significantly
lower than the average NCV error ($p_{0.5 m/s} = 0.0039$. Wilcoxon
one-sided rank-sign test) for the slowest speed ($0.5 m/s$). In other
words, assuming symmetry (CCV) produces a single step-to-step model
that has greater predictive power than is achieved by refining the
analysis to produce separate $\lbrace L \mapsto R \rbrace$ and
$\lbrace R \mapsto L \rbrace$ step maps.

\subsubsection{Stride Maps}

We analyzed the dynamical symmetry and the statistical consequences of
symmetric modeling on the stride-to-stride transitions. Similar to
before, we analyzed two different $\left( \mathcal{N} , \mathcal{M}
\right)$ combinations, $\left(\mathcal{N},\mathcal{M} \right) =
\left(\lbrace L \mapsto L \rbrace , \lbrace R \mapsto R
  \rbrace\right)$ and $\left( \mathcal{N},\mathcal{M} \right) =
\left(\lbrace R \mapsto R \rbrace , \lbrace L \mapsto L \rbrace
\right)$ (Table \ref{tab:cv}). And again, for each category of
cross validation, we averaged the CV errors for both combinations of
normal and mirrored data.  As in the previous section, we first
compared NCV and MCV errors to test if the stride-to-stride dataset is
statistically asymmetric. The NCV and CCV errors were also compared 
to contrast the symmetric and asymmetric modeling approaches.

\refig{fig:CrossCrossBoth}(b) compares the MCV and CCV errors to NCV
error for stride-to-stride data. MCV errors were higher (on average)
than NCV errors for all speeds, and these differences were
statistically significant ($p_{1.5 m/s} = 0.0391$, $p_{1 m/s} =
0.0117$ and $p_{0.5 m/s} = 0.0039$; paired one-sided Wilcoxon
rank-sign test). This shows that our dataset is dynamically asymmetric
between $L \mapsto L$ and $R \mapsto R$. 

However, the comparison of NCV and CCV errors in stride-to-stride dataset
is more striking than in the step-to-step case in that CCV errors were
statistically significantly lower than the NCV errors at all three
speeds ($p_{1.5 m/s} = 0.004$, $p_{1 m/s} = 0.012$ and $p_{0.5 m/s}
= 0.004$; paired one-sided Wilcoxon rank-sign test).

\subsubsection{Model Uncertainty}
\label{sec:consistency}

Cross-validation errors are powerful metrics for comparing the
effectiveness of symmetric and asymmetric modeling
approaches. However, if two models have similar CVEs, the next thing
to address is how well the data constrain the two models---i.e.\ how
much \emph{uncertainty} there is in the model
parameters~\cite{madhavclosed-loop2013}.  This was particularly important
for our step-to-step data because symmetric and asymmetric modeling
produced indistinguishable CVEs for $1.5$ and $1.0 m/s$ walking. This
implies that both modeling approaches are equally powerful from the
perspective of CVE. However, the parameters of the fitted section map
model may exhibit greater variability for the asymmetric modeling
approach since it uses less data for fitting.

In order to measure the uncertainty of the models, we adopted following
metric:
\begin{align}
\label{eq:consistency}
  \Xi = \sum\limits_{i=1}^d \sum\limits_{j=1}^d \sigma_{ij}^2 
\end{align}
where $\sigma_{ij}^2$ is the sample variance of $a_{ij}$, i.e. element
at the $i^{th}$ row and $j^{th}$ column of the section map
$A_\mathcal{F}$ fit during Monte-Carlo iterations of the extended CV
method. Symmetric model uncertainty was computed using the fitted matrix
samples of the CCV method. Model uncertainties of the $\lbrace L \mapsto R
\rbrace$ and $\lbrace R \mapsto L \rbrace$ (and
$\lbrace L \mapsto L \rbrace$ and $\lbrace R \mapsto R \rbrace$) maps
were averaged to have a single asymmetric model uncertainty for step
maps (and stride maps).

\begin{figure*}[t!] \centering
  \includegraphics[scale=1]{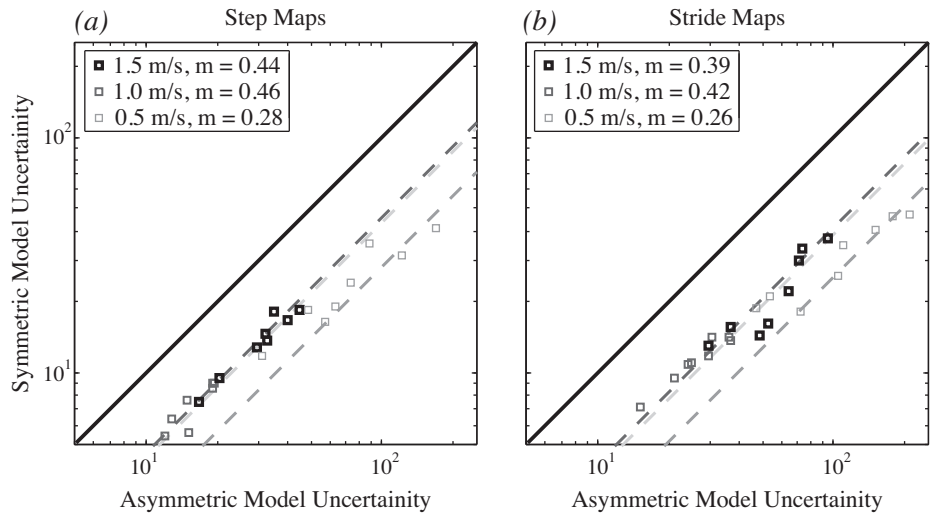}    
  \caption{Model uncertainty at all three speeds was lower when
    symmetry was assumed in both (A) step-to-step and (B)
    stride-to-stride maps. Each marker compares the model uncertainty with asymmetry
    and with symmetry of a single individual. Dashed line denotes the
    best fitted line (passing trough the origin) to the comparison
    markers. The percentage improvement is given by $(1-m) \times
    100$, where $m$ is the slope of the fitted line.}
  \label{fig:consistency}
\end{figure*}

We found that by neglecting asymmetry and fitting a single return map,
there was a substantial reduction in model uncertainty for both the
step-to-step and stride-to-stride data 
(Figure~\ref{fig:consistency}). Thus, even though in a few cases, the CV
errors were similar for NCV and CCV, the models produced using CCV
(that is, neglecting asymmetry and pooling the data) are substantially
less variable.

For step maps, assuming symmetry substantially lowers model
uncertainty; we saw $56\%$, $54\%$, and $72\%$ improvement with
symmetric approach for speeds 1.5, 1.0 and 0.5 m/s respectively. All
improvements were statistically significant ($p = 0.0039$, one-sided
Wilcoxon signed-rank test).  We observed the same trend with stride
maps: $61\%$, $58\%$, and $74\%$ improvement with symmetric approach
for speeds 1.5, 1.0 and 0.5 m/s respectively ($p = 0.0039$, one-sided
Wilcoxon signed-rank test). These results are illustrated in 
Figure~\ref{fig:consistency}.

\subsection{Step Return Maps vs. Stride Return Maps}

One of the advantages of assuming dynamical bilateral symmetry (i.e.\
neglecting asymmetry) is that one step becomes the fundamental period
of the system; the mapping from step to step defines the return map of
the dynamics. On the contrary, if we embrace the asymmetry, the stride
becomes the fundamental period. The disadvantage of using
stride-to-stride return maps compared to step-to-step maps is a potential 
loss of signal-to-noise ratio due to the fact that stride maps reduce the
temporal resolution. Thus, one can expect that stride-to-stride return maps 
would have lower predictive power in the CV setting.

In order to compare the predictive powers of step and stride return maps,
we analyzed the CVEs by assuming symmetry and fitting lumped return maps 
to both step and stride data. Specifically, we compared the CCV errors of step and
stride data in our method. In order to estimate CCV error of step/stride
return map, we took the mean of CCV errors of $\lbrace L \mapsto R
\rbrace$ and $\lbrace R \mapsto L \rbrace$/$\lbrace L \mapsto L \rbrace$ and $\lbrace R \mapsto R
\rbrace$.

\begin{figure}[ht] \centering
  \includegraphics[scale=1]{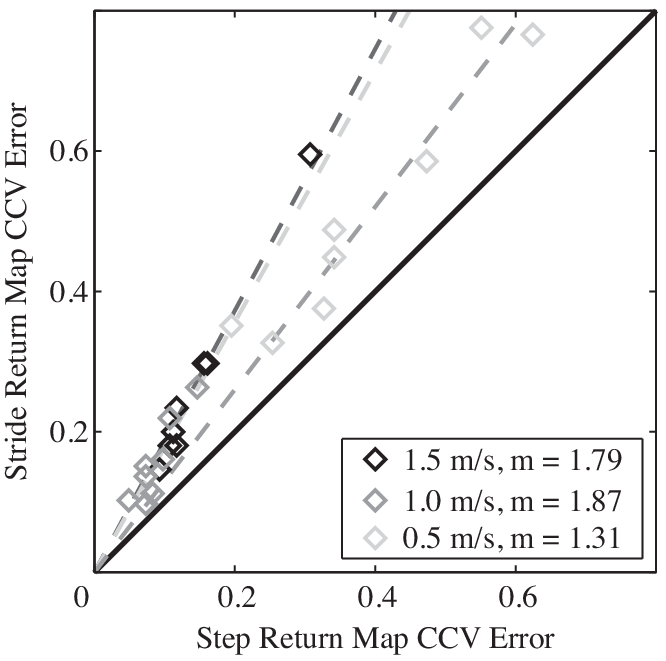}
  \caption{Illustration of the CCV errors of step return
    maps and stride return maps. Each marker compares the
    CCV errors of step and stride return maps of a single
    individual. Dashed lines illustrate the best fitted lines
    (passing trough the origin) to the comparison
    markers with $m$ being the associated slope of the line.
  }
  \label{fig:crossErrors}
\end{figure}

The results illustrated in \refig{fig:crossErrors} show that there is
a dramatic signal-to-noise ratio loss with stride-to-stride return
maps and that step-to-step return maps have more predictive power in the CV
setting. CCV errors of stride return maps are significantly higher
than the ones with step return maps: $79\%$, $87\%$, and $31\%$ more
CV error with stride return maps for speeds 1.5, 1.0, and 0.5 m/s
respectively. The differences are statistically significant ($p =
0.0039$, one-sided Wilcoxon signed-rank test).

\section{DISCUSSION}

In this paper, we focused our attention on bilateral dynamic asymmetry
in human walking.  Specifically, we introduced a statistical framework
based on applying cross validation techniques and fitting linear maps
to the data associated with the heel strike events. Our statistical
methods allowed us to examine the ``wrongness'' and ``usefulness'' of
neglecting bilateral dynamic asymmetry.

We applied our methods to data obtained from eight different
individuals walking at three different speeds. Based on the
results obtained with this data set, we observed that dynamical
asymmetry in walking is significant and statistically
distinguishable. These results underscore what several studies have
previously observed on steady-state
parameters~\cite{allardsimultaneous1996,lawmicrocomputer1987,stefanyshynright1994,damholtasymmetry1978,gundersenbilateral1989,ounpuubilateral1989,forczekevaluation2012,lathroplambachevidence2014}.

Despite the existence of significant asymmetry, we show that ignoring this and modeling 
human walking dynamics as symmetric produces significantly more consistent 
models (Figure~\ref{fig:consistency}). Moreover, the predictive power
of these symmetric models is higher than (or at worst equal to) their
asymmetric counterparts (Figure~\ref{fig:CrossCrossBoth}). This shows 
that neglecting bilateral asymmetry---an inescapable characteristic of the 
human form---not only provides modeling convenience but, more importantly, 
produces better models in terms of consistency and predictive power. It is 
not only ``OK'' to neglect asymmetry; in some cases, it is better. 

One should also note that the slight differences between two
``symmetrically'' placed sensors (e.g.\ load cells) can generate an
appearance of asymmetry that is not related to the actual
system. These asymmetries can affect both limit-cycle symmetry as well
as introduce dynamical asymmetries. Fortunately, despite possible
measurement asymmetries that would likely exacerbate asymmetries in
modeling, a symmetric dynamic model was still preferable for our
data. Another limitation to the present study is that participants
held their \textit{arms crossed} while they walked on a \textit{treadmill}, both of
which can affect gait \cite{ortegaeffects2008,dingwelllocal2001}.  As
instrumentation improves, the questions addressed in this paper can be
revisited in unconstrained and/or overground walking and it would be
interesting if walking were to be either more \emph{or} less symmetric
in those cases.

Even though we applied our methods to human walking data, they are
directly applicable to a wide range of rhythmic dynamical systems in
biology and robotics. Specifically, we are interested in behaviors that
exhibit alternating (out of phase) gait patterns but are symmetric via
reversing the left--right axis for half the stride. This class
includes bipedal walking, running, and sprinting \cite{mannbiomechanics1980,collinsefficient2005};
quadrupedal walking, trotting, and pacing
\cite{collinscoupled1993,buehlerscout1998}; hexapedal alternating
tripod gait \cite{fullmechanics1990,saranlirhex:2001}; and even
swimming \cite{sfakiotakisreview1999,ijspeertswimming2007}.

In the context of robotics, our methods can be used for diagnostics
and calibration since symmetry is considered a desirable property in
the design and development of robotic systems. Asymmetric robotic
gaits can potentially increase energy expenditure, reduce performance,
and introduce a steering bias, hindering the control and operation of
the robot. It may be possible to eliminate this steering bias by using
existing gait adaptation methods
\cite{weingartenautomated2004,gallowayexperimental2011} which, to
date, requires external instrumentation and specialized arenas.
However, our method relies on only internal kinematic measurements
which are directly available in most robotic systems, and so perhaps
the methods presented in this paper can be used to develop fast and
effective calibration methods for field robotics.

On the biological side, there is scientific value in investigating
dynamical symmetry across species. Models of biological locomotion can
be decomposed into two components: the mechanics of the locomotion
(plant), and the neural feedback (controller)~\cite{rothcomparative2014}. A ``less wrong'' model
of the plant provides better understanding of the controller, and
vice-versa~\cite{cowancritical2007,cowanfeedback2014}. The locomotor
pattern of a behaving animal is the closed-loop interaction of the
plant and controller. Investigating dynamical symmetry (or asymmetry)
in the locomotor gait as well as symmetry (or asymmetry) of the
kinematics allows us to better predict the structure of the
corresponding neural controller.

With regard to human health in particular, our tools may be useful for
understanding motor deficits during locomotion. Specifically, these
methods provide an important extension to those that center on
kinematic symmetry and its relations to human physiology
\cite{reismanlocomotor2007,finleylearning2013}. Individuals with
damage to the musculoskeletal system or nervous system often use
asymmetric kinematic walking patterns (e.g. amputees, stroke
patients). The kinematic asymmetry can be in the amount of time
standing on one leg versus the other, the extent of limb movements, or
some combination. An understanding of the underlying dynamical
asymmetry (or even symmetry) in these cases would provide more
information about the nature of the deficit, and perhaps suggest new
targets for focusing rehabilitation treatments.

Finally, an interesting extension of our methods would be analyzing
\emph{dynamical} asymmetry in gaits with categorically asymmetric
steady-state kinematics, such as quadrupedal galloping and
bounding. The steady-state limit-cycles of such gaits are obviously
asymmetric, but the dynamics around those limit-cycles may be
symmetric (enough).

\section{Data accessibility}

The data and codes used for drafting this paper will be able 
to be accessed via DOI:10.7281/T15Q4T12 (URL: \url{dx.doi.org/10.7281/T15Q4T12}).
The data used in this paper are the from the baseline walking portion 
in previous work by Long et al. \cite{longmarching-walking2015}.

\section{Authors' Contributions}

MMA developed the statistical tools, analyzed and interpreted the data, and drafted/revised the manuscript. 
SF, MSM, and NJC interpreted the data, and drafted/revised the manuscript. 
AL designed and performed the experiments, and revised the manuscript.
AJB designed the experiments, interpreted the data, and revised the manuscript.

\section{Acknowledgments and Funding Statement}

This material is based upon work supported by the National Science
Foundation (NSF) under grants 0845749 and 1230493 to N.~J.~Cowan and
by the National Institutes of Health (NIH) under grant R01-HD048741 to
A.~J.~Bastian. We also thank the reviewers for their insightful
comments and suggestions.

\clearpage

\section*{Tables}

\begin{table}[ht]
  \caption{Catalog of normal, $\mathcal{N}$, and associated mirrored, $\mathcal{M}$,
    datasets combinations used in our CV analysis.}
  \label{tab:cv}
  \centering
  \begin{tabular}{|c|c|c|} \hline
  & $\mathcal{N}$ & $\mathcal{M}$ 
  \\ \hline 
  Step &  $\begin{aligned}
      \lbrace L \mapsto R \rbrace
      \\ 
      \lbrace R \mapsto L \rbrace
    \end{aligned}$ & 
    $\begin{aligned}
      \lbrace R \mapsto L \rbrace
      \\ 
      \lbrace L \mapsto R \rbrace
    \end{aligned}$
  \\   \hline
  Stride &  $\begin{aligned}
      \lbrace L \mapsto L \rbrace
      \\ 
      \lbrace R \mapsto R \rbrace
    \end{aligned}$ & 
    $\begin{aligned}
      \lbrace R \mapsto R \rbrace
      \\ 
      \lbrace L \mapsto L \rbrace
    \end{aligned}$
  \\   \hline
  \end{tabular}
\end{table}

\section*{Figures}

\begin{thebibliography}{10}
\expandafter\ifx\csname urlstyle\endcsname\relax
  \providecommand{\doi}[1]{doi:\discretionary{}{}{}#1}\else
  \providecommand{\doi}{doi:\discretionary{}{}{}\begingroup
  \urlstyle{rm}\Url}\fi

\bibitem{finnertyorigins2004}
Finnerty, J.~R., Pang, K., Burton, P., Paulson, D. \& Martindale, M.~Q., 2004
  Origins of bilateral symmetry: Hox and dpp expression in a sea anemone.
\newblock \emph{Science} \textbf{304}, 1335--1337.
\newblock (\doi{10.1126/science.1091946}).

\bibitem{sadeghisymmetry2000}
Sadeghi, H., Allard, P., Prince, F. \& Labelle, H., 2000 {Symmetry and limb
  dominance in able-bodied gait: a review}.
\newblock \emph{Gait Posture} \textbf{12}, 34--45.
\newblock ISSN 09666362.
\newblock (\doi{10.1016/S0966-6362(00)00070-9}).

\bibitem{hsiaoreview2010}
Hsiao-Wecksler, E.~T., Polk, J.~D., Rosengren, K.~S., Sosnoff, J.~J. \& Hong,
  S., 2010 A review of new analytic techniques for quantifying symmetry in
  locomotion.
\newblock \emph{Symmetry} \textbf{2}, 1135--1155.
\newblock (\doi{10.3390/sym2021135}).

\bibitem{finleylearning2013}
Finley, J.~M., Bastian, A.~J. \& Gottschall, J.~S., 2013 Learning to be
  economical: the energy cost of walking tracks motor adaptation.
\newblock \emph{J Physiol} \textbf{591}, 1081--1095.
\newblock (\doi{10.1113/jphysiol.2012.245506}).

\bibitem{laigait2001}
Lai, K.-A., Lin, C.-J., Jou, I., Su, F.-C. \emph{et~al.}, 2001 Gait analysis
  after total hip arthroplasty with leg-length equalization in women with
  unilateral congenital complete dislocation of the hip--comparison with
  untreated patients.
\newblock \emph{J Orthop Res} \textbf{19}, 1147--1152.

\bibitem{matteswalking2000}
Mattes, S.~J., Martin, P.~E. \& Royer, T.~D., 2000 Walking symmetry and energy
  cost in persons with unilateral transtibial amputations: matching prosthetic
  and intact limb inertial properties.
\newblock \emph{Arch Phys Med Rehabil} \textbf{81}, 561--568.
\newblock (\doi{10.1016/S0003-9993(00)90035-2}).

\bibitem{hannahkinematic1984}
Hannah, R., Morrison, J. \& Chapman, A., 1984 Kinematic symmetry of the lower
  limbs.
\newblock \emph{Arch Phys Med Rehabil} \textbf{65}, 155--158.

\bibitem{forczekevaluation2012}
Forczek, W. \& Staszkiewicz, R., 2012 {An evaluation of symmetry in the lower
  limb joints during the able-bodied gait of women and men.}
\newblock \emph{J Hum Kinet} \textbf{35}, 47--57.
\newblock ISSN 1640-5544.
\newblock (\doi{10.2478/v10078-012-0078-5}).

\bibitem{karamanidissymmetry2003}
Karamanidis, K., Arampatzis, A. \& Bruggemann, G.~P., 2003 Symmetry and
  reproducibility of kinematic parameters during various running techniques.
\newblock \emph{Med Sci Sports Exerc} \textbf{35}, 1009--1016.
\newblock (\doi{10.1249/01.MSS.0000069337.49567.F0}).

\bibitem{reismanlocomotor2007}
Reisman, D.~S., Wityk, R., Silver, K. \& Bastian, A.~J., 2007 Locomotor
  adaptation on a split-belt treadmill can improve walking symmetry
  post-stroke.
\newblock \emph{Brain} \textbf{130}, 1861--1872.
\newblock (\doi{10.1093/brain/awm035}).

\bibitem{stefanyshynright1994}
Stefanyshyn, D.~J. \& Engsberg, J.~R., 1994 Right to left differences in the
  ankle joint complex range of motion.
\newblock \emph{Med Sci Sports Exerc} \textbf{26}, 551--555.

\bibitem{gundersenbilateral1989}
Gundersen, L.~A., Valle, D.~R., Barr, A.~E., Danoff, J.~V., Stanhope, S.~J. \&
  Snyder-Mackler, L., 1989 {Bilateral analysis of the knee and ankle during
  gait: an examination of the relationship between lateral dominance and
  symmetry.}
\newblock \emph{Phys Ther} \textbf{69}, 640--50.
\newblock ISSN 0031-9023.

\bibitem{lawmicrocomputer1987}
Law, H., 1987 Microcomputer-based low-cost method for measurement of spatial
  and temporal parameters of gait.
\newblock \emph{J Biomed Eng} \textbf{9}, 115--120.

\bibitem{choderafootprint1973}
Chodera, J. \& Levell, R., 1973 \emph{Footprint patterns during walking}, pp.
  81--90.
\newblock Perspectives in Biomedical Engineering. Baltimore: University Park
  Press.

\bibitem{choderaanalysis1974}
Chodera, J., 1974 Analysis of gait from footprints.
\newblock \emph{Physiotherapy} \textbf{60}, 179.

\bibitem{hamillground1984}
Hamill, J., Bates, B. \& Knutzen, K., 1984 Ground reaction force symmetry
  during walking and running.
\newblock \emph{Res Q Exerc Sport} \textbf{55}, 289--293.

\bibitem{menardcomparative1992}
Menard, M.~R., McBride, M.~E., Sanderson, D.~J. \& Murray, D.~D., 1992
  Comparative biomechanical analysis of energy-storing prosthetic feet.
\newblock \emph{Arch Phys Med Rehabil} \textbf{73}, 451--58.

\bibitem{vanderstraatensymmetry1978}
Van~der Straaten, J. \& Scholton, P., 1978 Symmetry and periodicity in gait
  patterns of normal and hemiplegic children.
\newblock In \emph{Biomechanics VI: proceedings of the Sixth International
  Congress of Biomechanics, Copenhagen, Denmark}, volume~1, p. 287. University
  Park Press.

\bibitem{herzogasymmetries1989}
Herzog, W., Nigg, B.~M., Read, L.~J. \& Olsson, E., 1989 {Asymmetries in ground
  reaction force patterns in normal human gait.}
\newblock \emph{Med Sci Sports Exerc} \textbf{21}, 110--4.
\newblock ISSN 0195-9131.

\bibitem{carlsookinetic1973}
Carls{\"o}{\"o}, S., Dahl{\"o}f, A. \& Holm, J., 1973 Kinetic analysis of the
  gait in patients with hemiparesis and in patients with intermittent
  claudication.
\newblock \emph{Scand J Rehabil Med} \textbf{6}, 166--179.

\bibitem{arsenaultthere1986}
Arsenault, A., Winter, D. \& Marteniuk, R., 1986 Is there a `normal' profile of
  emg activity in gait?
\newblock \emph{Med Biol Eng Comput} \textbf{24}, 337--343.

\bibitem{marksanalysis1958}
Marks, M. \& Hirschberg, G.~G., 1958 Analysis of the hemiplegic gait.
\newblock \emph{Ann NY Acad Sci} \textbf{74}, 59--77.

\bibitem{ounpuubilateral1989}
{\~O}unpuu, S. \& Winter, D.~A., 1989 Bilateral electromyographical analysis of
  the lower limbs during walking in normal adults.
\newblock \emph{Electroencephalogr Clin Neurophysiol} \textbf{72}, 429--438.

\bibitem{damholtasymmetry1978}
Damholt, V. \& Termansen, N., 1978 Asymmetry of plantar flexion strength in the
  foot.
\newblock \emph{Acta Orthop} \textbf{49}, 215--219.

\bibitem{balakrishanintegral1982}
Balakrishan, S. \& Thornton-Trump, A., 1982 Integral parameters in human
  locomotion.
\newblock In \emph{Proceeding of the Second Biannual Conference of the Canadian
  Society for Biomechanics, Human Locomotion II}, pp. 12--3.

\bibitem{vaughanjoint1996}
Vaughan, C.~L., 1996 {Are joint torques the Holy Grail of human gait analysis?}
\newblock \emph{Hum Movement Sci} \textbf{15}, 423--443.

\bibitem{lathroplambachevidence2014}
Lathrop-Lambach, R.~L., Asay, J.~L., Jamison, S.~T., Pan, X., Schmitt, L.~C.,
  Blazek, K., Siston, R.~a., Andriacchi, T.~P. \& Chaudhari, A. M.~W., 2014
  {Evidence for joint moment asymmetry in healthy populations during gait.}
\newblock \emph{Gait Posture} \textbf{40}, 526--31.
\newblock ISSN 1879-2219.
\newblock (\doi{10.1016/j.gaitpost.2014.06.010}).

\bibitem{crowecharacterization1995}
Crowe, A., Schiereck, P., De~Boer, R. \& Keessen, W., 1995 Characterization of
  human gait by means of body center of mass oscillations derived from ground
  reaction forces.
\newblock \emph{IEEE Trans Biomed Eng} \textbf{42}, 293--303.

\bibitem{crowecharacterization1993}
Crowe, A., Schiereck, P., de~Boer, R. \& Keessen, W., 1993 Characterization of
  gait of young adult females by means of body centre of mass oscillations
  derived from ground reaction forces.
\newblock \emph{Gait Posture} \textbf{1}, 61--68.

\bibitem{giakastime1997}
Giakas, G. \& Baltzopoulos, V., 1997 Time and frequency domain analysis of
  ground reaction forces during walking: an investigation of variability and
  symmetry.
\newblock \emph{Gait Posture} \textbf{5}, 189--197.
\newblock (\doi{10.1016/S0966-6362(96)01083-1}).

\bibitem{boxempirical1987}
Box, G. E.~P. \& Draper, N.~R., 1987 \emph{Empirical Model-Building and
  Response Surfaces}.
\newblock Wiley.

\bibitem{altendorferstability2004a}
Altendorfer, R., Koditschek, D.~E. \& Holmes, P., 2004 Stability analysis of
  legged locomotion models by symmetry-factored return maps.
\newblock \emph{Int J Robot Res} \textbf{23}, 979--999.
\newblock (\doi{10.1177/0278364904047389}).

\bibitem{chevallereauasymptotically2009}
Chevallereau, C., Grizzle, J.~W. \& Shih, C.-L., 2009 Asymptotically stable
  walking of a five-link underactuated 3-d bipedal robot.
\newblock \emph{IEEE Trans Robot} \textbf{25}, 37--50.
\newblock (\doi{10.1109/TRO.2008.2010366}).

\bibitem{ankarali_saranli.chaos2010}
Ankarali, M.~M. \& Saranli, U., 2010 Stride-to-stride energy regulation for
  robust self-stability of a torque-actuated dissipative spring-mass hopper.
\newblock \emph{Chaos} \textbf{20}.
\newblock (\doi{10.1063/1.3486803}).

\bibitem{depenn2015}
De, A. \& Koditschek, D.~E., 2015 The penn jerboa: A platform for exploring
  parallel composition of templates.
\newblock \emph{arXiv preprint arXiv:1502.05347} .

\bibitem{hurmuzlumeasurement1994}
Hurmuzlu, Y. \& Basdogan, C., 1994 On the measurement of dynamic stability of
  human locomotion.
\newblock \emph{Trans ASME} \textbf{116}, 30--36.

\bibitem{seyfarthswing-leg2003}
Seyfarth, A., Geyer, H. \& Herr, H., 2003 Swing-leg retraction: a simple
  control model for stable running.
\newblock \emph{J Exp Biol} \textbf{206}, 2547--2555.
\newblock (\doi{10.1242/​jeb.00463}).

\bibitem{holmesdynamics2006}
Holmes, P.~J., Full, R.~J., Koditschek, D.~E. \& Guckenheimer, J., 2006 The
  dynamics of legged locomotion: Models, analyses, and challenges.
\newblock \emph{SIAM Rev} \textbf{48}, 207--304.
\newblock (\doi{10.1137/S0036144504445133}).

\bibitem{revzenfinding2012}
Revzen, S. \& Guckenheimer, J.~M., 2011 Finding the dimension of slow dynamics
  in a rhythmic system.
\newblock \emph{J R Soc Interface} \textbf{9}, 957--971.
\newblock (\doi{10.1098/rsif.2011.0431}).

\bibitem{ankaralihaptic2014}
Ankarali, M.~M., {\c S}en, H.~T., De, A., Okamura, A.~M. \& Cowan, N.~J., 2014
  Haptic feedback enhances rhythmic motor control by reducing variability, not
  improving convergence rate.
\newblock \emph{J Neurophysiol} \textbf{111}, 1286--1299.
\newblock (\doi{10.1152/jn.00140.2013}).

\bibitem{guckenheimernonlinear1991}
Guckenheimer, J. \& Holmes, P., 1991 \emph{Nonlinear Oscillations, Dynamical
  Systems, and Bifurcations of Vector Fields}.
\newblock Springer.

\bibitem{muircomplete1999}
Muir, G.~D. \& Whishaw, I.~Q., 1999 Complete locomotor recovery following
  corticospinal tract lesions: measurement of ground reaction forces during
  overground locomotion in rats.
\newblock \emph{Behav Brain Res} \textbf{103}, 45--53.
\newblock (\doi{10.1016/S0166-4328(99)00018-2}).

\bibitem{pourcelotkinematic1997}
Pourcelot, P., Audigie, F., Degueurce, C., Denoix, J. \& Geiger, D., 1997
  Kinematic symmetry index: a method for quantifying the horse locomotion
  symmetry using kinematic data.
\newblock \emph{Vet Res} \textbf{28}, 525.

\bibitem{millercontinuous2008}
Miller, R.~H., Meardon, S.~A., Derrick, T.~R. \& Gillette, J.~C., 2008
  Continuous relative phase variability during an exhaustive run in runners
  with a history of iliotibial band syndrome.
\newblock \emph{J Appl Biomech} \textbf{24}, 262--270.

\bibitem{kurzdifferences2012}
Kurz, M.~J., Arpin, D.~J. \& Corr, B., 2012 Differences in the dynamic gait
  stability of children with cerebral palsy and typically developing children.
\newblock \emph{Gait Posture} \textbf{36}, 600--604.

\bibitem{donelanmechanical2002}
Donelan, J.~M., Kram, R. \& Kuo, A.~D., 2002 Mechanical work for step-to-step
  transitions is a major determinant of the metabolic cost of human walking.
\newblock \emph{J Exp Biol} \textbf{205}, 3717--3727.

\bibitem{ankarali_saranli_AR2011}
Ankarali, M. \& Saranli, U., 2011 Control of underactuated planar pronking
  through an embedded spring-mass hopper template.
\newblock \emph{Auton Robot} \textbf{30}, 217--231.
\newblock ISSN 0929-5593.
\newblock (\doi{10.1007/s10514-010-9216-x}).

\bibitem{saranli_arslan_ankarali_morgul.nd2010}
Saranli, U., Arslan, O., Ankarali, M.~M. \& Morgul, O., 2010 Approximate
  analytic solutions to non-symmetric stance trajectories of the passive
  spring-loaded inverted pendulum with damping.
\newblock \emph{Nonlinear Dynam} \textbf{62}, 729--742.
\newblock (\doi{10.1007/s11071-010-9757-8}).

\bibitem{blickhansimilarity1993}
Blickhan, R. \& Full, R.~J., 1993 Similarity in multilegged locomotion:
  Bouncing like a monopode.
\newblock \emph{J Comp Physiol A} \textbf{173}, 509--517.
\newblock (\doi{10.1007/BF00197760}).

\bibitem{westervelthybrid2003}
Westervelt, E., Grizzle, J. \& Koditschek, D., 2003 Hybrid zero dynamics of
  planar biped walkers.
\newblock \emph{IEEE Trans Autom Control} \textbf{48}, 42--56.
\newblock (\doi{10.1109/TAC.2002.806653}).

\bibitem{duindammodeling2009}
Duindam, V. \& Stramigioli, S., 2009 Modeling and analysis of walking robots.
\newblock In \emph{Modeling and Control for Efficient Bipedal Walking Robots},
  pp. 93--127. Springer.

\bibitem{leetemplates2008}
Lee, J., Sponberg, S.~N., Loh, O.~Y., Lamperski, A.~G., Full, R.~J. \& Cowan,
  N.~J., 2008 Templates and anchors for antenna-based wall following in
  cockroaches and robots.
\newblock \emph{IEEE Trans Robot} \textbf{24}, 130--143.
\newblock (\doi{10.1109/TRO.2007.913981}).

\bibitem{ankaralivariability2015}
Ankarali, M.~M., 2015 \emph{Variability, Symmetry, and Dynamics in Human
  Rhythmic Motor Control}.
\newblock Ph.D. thesis, Johns Hopkins University.

\bibitem{shaolinear1993}
Shao, J., 1993 Linear model selection by cross-validation.
\newblock \emph{J Am Statist Assoc} \textbf{88}, 486--494.

\bibitem{madhavclosed-loop2013}
Madhav, M.~S., Stamper, S.~A., Fortune, E.~S. \& Cowan, N.~J., 2013 Closed-loop
  stabilization of the jamming avoidance response reveals its locally unstable
  and globally nonlinear dynamics.
\newblock \emph{J Exp Biol} \textbf{216}, 4272--4284.
\newblock (\doi{10.1242/jeb.088922}).

\bibitem{raolinear2005}
Rao, C. \& Wu, Y., 2005 Linear model selection by cross-validation.
\newblock \emph{J Stat Plan Inference} \textbf{128}, 231--240.
\newblock (\doi{10.1016/j.jspi.2003.10.004}).

\bibitem{yangconsistency2007}
Yang, Y., 2007 Consistency of cross validation for comparing regression
  procedures.
\newblock \emph{Ann Stat} pp. 2450--2473.

\bibitem{arlotsurvey2010}
Arlot, S. \& Celisse, A., 2010 A survey of cross-validation procedures for
  model selection.
\newblock \emph{Stat Surv} \textbf{4}, 40--79.

\bibitem{collinscoupled1993}
Collins, J.~J. \& Stewart, I.~N., 1993 Coupled nonlinear oscillators and the
  symmetries of animal gaits.
\newblock \emph{J Nonlinear Sci} \textbf{3}, 349--392.
\newblock (\doi{10.1007/BF02429870}).

\bibitem{allardsimultaneous1996}
Allard, P., Lachance, R., Aissaoui, R. \& Duhaime, M., 1996 Simultaneous
  bilateral 3-d able-bodied gait.
\newblock \emph{Hum Movement Sci} \textbf{15}, 327--346.

\bibitem{ortegaeffects2008}
Ortega, J.~D., Fehlman, L.~A. \& Farley, C.~T., 2008 Effects of aging and arm
  swing on the metabolic cost of stability in human walking.
\newblock \emph{J Biomech} \textbf{41}, 3303--3308.

\bibitem{dingwelllocal2001}
Dingwell, J., Cusumano, J., Cavanagh, P. \& Sternad, D., 2001 Local dynamic
  stability versus kinematic variability of continuous overground and treadmill
  walking.
\newblock \emph{J Biomech Eng} \textbf{123}, 27--32.

\bibitem{mannbiomechanics1980}
Mann, R.~A. \& Hagy, J., 1980 Biomechanics of walking, running, and sprinting.
\newblock \emph{Am J Sports Med} \textbf{8}, 345--350.
\newblock (\doi{10.1177/036354658000800510}).

\bibitem{collinsefficient2005}
Collins, S., Ruina, A., Tedrake, R. \& Wisse, M., 2005 Efficient bipedal robots
  based on passive-dynamic walkers.
\newblock \emph{Science} \textbf{307}, 1082--5.
\newblock (\doi{10.1126/science.1107799}).

\bibitem{buehlerscout1998}
Buehler, M., Battaglia, R., Cocosco, A., Hawker, G., Sarkis, J. \& Yamazaki,
  K., 1998 Scout: A simple quadruped that walks, climbs, and runs.
\newblock In \emph{Proc IEEE Int Conf Robot Autom}, volume~2, pp. 1707--1712.
  IEEE.

\bibitem{fullmechanics1990}
Full, R.~J. \& Tu, M.~S., 1990 Mechanics of six-legged runners.
\newblock \emph{J Exp Biol} \textbf{148}, 129--146.

\bibitem{saranlirhex:2001}
Saranli, U., Buehler, M. \& Koditschek, D.~E., 2001 {RH}ex: A simple and highly
  mobile hexapod robot.
\newblock \emph{Int J Robot Res} \textbf{20}, 616--631.
\newblock (\doi{10.1177/02783640122067570}).

\bibitem{sfakiotakisreview1999}
Sfakiotakis, M., Lane, D.~M. \& Davies, J. B.~C., 1999 Review of fish swimming
  modes for aquatic locomotion.
\newblock \emph{IEEE J Ocean Eng} \textbf{24}, 237--252.
\newblock (\doi{10.1109/48.757275}).

\bibitem{ijspeertswimming2007}
Ijspeert, A.~J., Crespi, A., Ryczko, D. \& Cabelguen, J.-M., 2007 From swimming
  to walking with a salamander robot driven by a spinal cord model.
\newblock \emph{Science} \textbf{315}, 1416--1420.
\newblock (\doi{10.1126/science.1138353}).

\bibitem{weingartenautomated2004}
Weingarten, J., Lopes, G. A.~D., Buehler, M., Groff, R.~E. \& Koditschek, D.,
  2004 Automated gait adaptation for legged robots.
\newblock In \emph{Proc IEEE Int Conf Robot Autom}, volume~3, pp. 2153--2158.
\newblock (\doi{10.1109/ROBOT.2004.1307381}).

\bibitem{gallowayexperimental2011}
Galloway, K., Clark, J., Yim, M. \& Koditschek, D., 2011 Experimental
  investigations into the role of passive variable compliant legs for dynamic
  robotic locomotion.
\newblock In \emph{Proc IEEE Int Conf Robot Autom}, pp. 1243--1249.
\newblock (\doi{10.1109/ICRA.2011.5979941}).

\bibitem{rothcomparative2014}
Roth, E., Sponberg, S. \& Cowan, N.~J., 2014 A comparative approach to
  closed-loop computation.
\newblock \emph{Curr Opin Neurobiol} \textbf{25}, 54--62.
\newblock (\doi{10.1016/j.conb.2013.11.005}).

\bibitem{cowancritical2007}
Cowan, N.~J. \& Fortune, E.~S., 2007 The critical role of locomotion mechanics
  in decoding sensory systems.
\newblock \emph{J Neurosci} \textbf{27}, 1123--1128.
\newblock (\doi{10.1523/JNEUROSCI.4198-06.2007}).

\bibitem{cowanfeedback2014}
Cowan, N.~J., Ankarali, M.~M., Dyhr, J.~P., Madhav, M.~S., Roth, E., Sefati,
  S., Sponberg, S., Stamper, S.~A., Fortune, E.~S. \& Daniel, T.~L., 2014
  Feedback control as a framework for understanding tradeoffs in biology.
\newblock \emph{Integr Comp Biol} \textbf{54}, 223--237.
\newblock ISSN 1540-7063.
\newblock (\doi{10.1093/icb/icu050}).

\bibitem{longmarching-walking2015}
Long, A.~W., Finley, J.~M. \& Bastian, A.~J., 2015 A Marching-Walking Hybrid Induces Step 
Length Adaptation and Transfers to Natural Walking.
\newblock \emph{J Physiol}.
\newblock ISSN 0022-3077.
\newblock (\doi{10.1152/jn.00779.2014}).

\end{thebibliography}
\end{document}